\documentclass{amsart}
\usepackage{graphicx}
\usepackage{psfrag, subfigure}

\newtheorem{theorem}{Theorem}[section]
\newtheorem{lemma}[theorem]{Lemma}

\theoremstyle{definition}
\newtheorem{definition}[theorem]{Definition}

\newtheorem{rk}{Remark}
\theoremstyle{remark}

\theoremstyle{corollary}
\newtheorem{corollary}[theorem]{Corollary}

\newcommand{\N}{\mbox{$\mathbb{N}$}}
\newcommand{\Z}{\mbox{$Z\!\!\!Z$}}

\newcommand{\R}{\mbox{$I\!\!R$}}

\numberwithin{equation}{section}

\newcommand{\dem}{{\textbf{Proof:\ }}}

\begin{document}

\title{About $C^{1}$-minimality of the hyperbolic Cantor sets.}
\author{Liane Bordignon, Jorge Iglesias and Aldo Portela}

\begin{abstract}
 In this work we prove that a $C^{1+\alpha}$-hyperbolic Cantor set contained in $S^1$, close to an affine Cantor set is not $C^{1}$-minimal.
\end{abstract}

\maketitle

\section{Introduction}
If $f:S^1 \rightarrow S^1$ is a $C^1$-diffeomorphism without
periodic points, there exists a unique set $\Omega (f) \subset
S^1$ minimal for $f$ (we say that $\Omega (f) $ is $C^1$-minimal
for $f$). In this case $\Omega (f)$ is either a Cantor set or
$S^1$. The $C^1$-minimal Cantor sets that are known are the
Denjoy's  examples and their conjugates (see \cite{D}). Examples of Cantor sets that are not $C^1$-minimal, are also known, see  \cite{mc}, \cite{IP}, \cite{N}, \cite{N2} and \cite{P}.

Let $I_1,...,I_k$, $k\geq 2$ be pairwise disjoint compact intervals in $\R$,
and let $L$ be a compact interval containing their union $I\equiv  I_1\cup \cdots\cup I_k$. Define
${\mathcal S}^{r}(I_1,...,I_k,L)$, $r\geq 1$, to be the set of $C^{r}$ functions $S:I\to L$ such that for $j=1,...,k$, $S(I_j)=L$ and $|S^{'}|>1$.
For $S\in  {\mathcal S}^{r}(I_1,...,I_k,L)$ define $$C_S=\{x\in I:\ S^{k}(x)\in I\mbox{ for all } k\in \Z^{+}\}.$$
These sets are Cantor sets, and are called $C^{r}$-\textit{hyperbolic}. If $S^{'}$ is locally constant, $C_S$ is called  \textit{affine} and if $S^{'}$ is globally  constant, $C_S$ is called \textit{linear}.
Let $S^{1}=\R/\Z$ and suppose $L\subset [0,1)$. We may assume that one of the complementary intervals of $C_S$ is $(1/2,1)$ and $L=[0,1/2]$.

In \cite{N}, A. N. Kercheval assert that the $C_S$ Cantor sets, $C^2$ close to an affine Cantor set, are not $C^1$-minimal, for $S$ increasing or decreasing and any $k$. However, in this work, only was proved it for $k=2$ and $S$ locally increasing. For $k>2$ or $k=2$ but $S$ not locally increasing it is not possible to generalize his proof.
In this work we will prove the statements given by Kercheval in the general context. Moreover we give a generalization for the case that the cantor sets $C_S$ are $C^{1+\alpha}$ close to a affine Cantor set.



We follow the notation of \cite{N}. If $S\in \mathcal{S}^{r}(I_1,...,I_k,L)$, let $\phi_i: L\to I_i $ , for all $1\leq i\leq k$,  be the $k$ branches of the inverse of $S$. For any choice $(i_1,...,i_n)\in \{1,\dots, k\}^n$
we write $\phi (i_1,...,i_n)$ for the composition $\phi _{i_{1}}\circ \cdots\circ \phi _{i_{n}}$. Call the interval $\phi (i_1,...,i_n)(L)$ an $n$-{\textit {block}} of $C_S$.  If  $T=\phi (i_1,...,i_n)(L)$,  $T(j_1,...,j_m)$ denote the interval
$$\phi (i_1,...,i_n,j_1,...,j_m)(L)\subset T.$$

If necessary, we will denote such blocks by $T^S$ and the map $\phi (i_1,...,i_n)$ by $\phi^S (i_1,...,i_n)$.

Call $J_i$, $1\leq i\leq k-1$, the interval between $I_i$ and $I_{i+1}$. For any $m\geq 0$,  denote $J_i(m)$ an interval $J$ such that $S^m(J)= J_i$ for some $ 1\leq i\leq k-1 $. The intervals $J_i(m)$ are called \textit{gaps} of \textit{level} $m$. Remark that $J_i(0)=J_i$.

Note that for each $m$ there exist $k^{m}$ gaps of level $m$. These gaps are the intervals
$\phi (i_1,...,i_m)(J_i)$ and we will  denote them by  $J_i(i_1,...,i_m)$ or $J^S_i(i_1,...,i_m)$, if necessary.

For $S\in {\mathcal S}^{r}(I_1,...,I_k,L)$, define the {\em nonlinearity} of $S$ to be $${\mathcal N}(S)=\underset{j=1,...,k}{\mbox{max}} \ {\mbox{}}\underset{x,y\in I_{j}}{\mbox{sup}}\mbox{log}\frac{S'(x)}{S'(y)}\cdot$$

For $M>0,\varepsilon >0$ and $\alpha >0$ let $${\mathcal A}(M,\varepsilon,\alpha )=\{S\in {\mathcal S}^{1+\alpha}(I_1,...,I_k,L): \ \  {\mathcal N}(S)<\varepsilon \mbox{ and } \ |\log \frac{|S'(x)|}{ |S'(y)|}|\leq M|x-y|^{\alpha}   \}.$$

Note that if $\varepsilon\to 0$ then $S\to R$ in the $C^{1}$ topology, where $R$ is the affine map of ${\mathcal S}^{r}(I_1,...,I_k,L)$ with $R'(x)S'(x) >0$ for all $x\in \cup I_i$.

In this work, we prove the following result:

\begin{theorem}\label{teorema}
Let $I_1,...,I_k,L$ be compact intervals in $[0,1)$ as above. Then for each $M>0$ and $\alpha >0$ there exists $\varepsilon >0$ such that if $S\in {\mathcal A}(M,\varepsilon,\alpha )$ then $C_S$ is not $C^{1}$-minimal.
\end{theorem}

\section{Proof of Theorem \ref{teorema}}

All along this section, $I_1,...,I_k,L$ will denote the intervals in $[0,1)$ as in the hypothesis of Theorem \ref{teorema}.

\begin{lemma}\label{lema1} Let  $R, S\in \mathcal{S}^1(I_1,...,I_k,L)$ such that $R'(x)S'(x) >0$ for all $x\in \cup I_i$. There exists a homeomorphism $\psi:L\to L$ such that $\psi(C_S)= C_R$ and $$\psi(\phi^S(i_1, \dots, i_m)(L))=\phi^R(i_1, \dots, i_m)(L),$$
for all $m\in \mathbb{N}$.
\end{lemma}
\dem  If $x\in J_i$ for some $1\leq i \leq k-1$ we define $\psi(x)=x$ and if $x\in \phi^S(i_1,\dots, i_n)(J_i)$, we define $\psi(x)= \phi^R(i_1,\dots, i_n)( S^n (x))$. Then we  extend  this continuous and monotonic function  to a homeomorphism $\psi :L\to L$ in the natural way.\hfill $\square$

\begin{definition} \label{gt} Let $S \in \mathcal{S}^{r}(I_1,...,I_k,L)$ and $f:S^1\to S^1$ be a $C^1$ diffeomorphism with $C_S$ as its minimal set.  If $T$ is an $m$-block, $T=\phi^S(i_1, \dots, i_m)(L)$, define
 \begin{eqnarray} u_T &=& \min \{ n\in \mathbb{N}\mid S^n (f(J_i (r)))= J_j, \ 1\leq i, j\leq k-1, J_i (r)\subset T \} \nonumber\\
&=&\max \{ n\in \mathbb{N}\mid S^{n-1} (f(T)) \subset \cup I_i \}.\nonumber
\end{eqnarray}

Also, define
\begin{equation}\label{defg}
g_{_T}= S^{u_T}\circ f_{\mid_{T}}.
\end{equation}

\end{definition}

The map defined in (\ref{defg}) can be increasing or decreasing, depending on $T$. Also, if $T$ is  an $m$-block, there exists $J_i (m+p) \subset T$ such that $ g_{_T}(J_i (m+p))= J_j$ for some $j \in \{1,..., k-1 \} $.

\begin{lemma}\label{lema2}
For all $\eta>0$, $\alpha>0$ and  $M>0$, there exists $\varepsilon >0$ such that if
$S\in {\mathcal A}(M,\varepsilon,\alpha )$  and $C_S$ is a minimal Cantor
set for some $C^1$ diffeomorphism $f_{_{S}}:S^1\to S^1$, there exists $m_0=m_0(f_{_{S}}) \in
\mathbb{N}$ such that if $T$ is  an $m$-block of $C_S$ with $m>m_0$, then
$$\mathcal{N}({g_{_T}}) <\eta,$$
where $g_{_T}$ is as in Definition \ref{gt}.
\end{lemma}

\dem Suppose  $S\in {\mathcal A}(M,\varepsilon_1,\alpha )$ for  fixed $I_1,\dots, I_k$ and  some positive numbers $M$, $\alpha$ and $\varepsilon_1$. Also, suppose $\varepsilon_1$ small enough that it assures  the existence of a $\sigma \in (0,1)$ such that
\begin{equation}\label{sigma}
|S'(x)|>\displaystyle\frac{1}{\sigma} ,\ \ \forall x\in I_1\cup \cdots \cup I_k \ {\mathrm{and}} \
\forall  S\in {\mathcal A}(M,\varepsilon_1,\alpha ).
\end{equation}

Then, by {\cite[Bounded Distortion Lemma]{N}}, for every $\delta >0$ there is $N\in \mathbb{N}$  such that if $m>N$ and $I$ is contained in an  $m$-block of $C_S$, \begin{equation}\label{delta}
\mathcal{N}({S^{m-N+1}}_{\mid_I}) <\delta ,  \ \  \forall  S\in {\mathcal A}(M,\varepsilon_1,\alpha ).
\end{equation}

Remark that if we replace only $\varepsilon_1$ for a smaller one, the inequality (\ref{sigma}) remains true for the same $\sigma$ as well as  the inequality (\ref{delta}) holds for the same $N$.

By {\cite[Lemma 1]{N}}, for any  $m$-block $T$,

 \begin{equation}\label{distorcao}  \mathcal{N}({g_{_T}})=\mathcal{N}({ S^{u_T}\circ f}_{\mid_{T}})
 \leq  \mathcal{N}({S^{u_T}}_{\mid_{f(T)}}) +\mathcal{N}(f_{\mid_{T}}) .
 \end{equation}

 Since $\mathcal{N}(S) <\varepsilon_1$,  by {\cite[Lemma 1]{N}}, if $0\leq N\leq u_T$,

 $$\mathcal{N}({S^{u_T}}_{\mid_{f(T)}})  \leq \mathcal{N}({S^{u_T-N}}_{\mid_{f(T)}})+\mathcal{N}({S^{N}}_{\mid_{S^{(u_{T}-N)}( f(T))}}).$$

Suppose $\delta>0 $ and $N$ as in the inequality (\ref{delta}). Remark that if $T$ is an $n$-block , $n$ great enough, then $f(T)$ is contained in an $m$-block with $m>N$. With these conditions, by inequality (\ref{delta}),
$$\mathcal{N}({S^{u_T-N}}_{\mid_{f(T)}})<\delta.$$
Also, since $\mathcal{N}(S) <\varepsilon_1$,  $\mathcal{N}({S^{N}}_{\mid_{S^{(u_T-N)}( f(T))}})\leq N\varepsilon_1$. Therefore
\begin{equation}\label{tt}
 \mathcal{N}({S^{u_T}}_{\mid_{f(T)}})\leq N\varepsilon_1 +\delta .
\end{equation}

  Let $\eta>0$ and take $\delta = \frac{\eta}{3}$. With this, fix $N$ satisfying the inequality (\ref{delta}).

 Now, let $\varepsilon>0$, $\varepsilon < \min \{ \varepsilon_1,\frac{\eta}{3N} \}$. Remark that we can replace $\varepsilon_1$ by $\varepsilon$ all along this proof, until here, and the inequalities remains true.

 Also, we take $m_0 \in \mathbb{N}$, $m_0=m_0(S,f)$,  $m_0>N$, such that if $T$ is an $m$-block with $m>m_0$,   then

 \begin{equation}\label{1111}
 \mathcal{N}(f_{\mid_{T}})< \frac{\eta}{3}.
\end{equation}
It is possible because  $f$ is a $C^1$ diffeomorphism.

 With these $N$, $\delta$, $\varepsilon$ and $m_0$,  if $T$ is an $m$-block with $m>m_0$, by inequalities (\ref{distorcao}), and (\ref{1111}),
 $$\mathcal{N}({g_{_T}})\leq N\varepsilon +\delta +\frac{\eta}{3} <\eta.$$
 \hfill $\square$

\begin{rk}\label{rk1}
For a given $\eta >0$, by the inequality (\ref{delta}) and {\cite[Lemma 1]{N}}, there exists $\varepsilon >0$ such that if $S \in {\mathcal A}(M,\varepsilon,\alpha )$ and $I$ is contained in an $m$-block then $$\mathcal{N}(S^{m}_{\mid_{ I}})<\eta .$$
\end{rk}

\begin{lemma}\label{lema33}
Let $M>0$ and $\alpha >0$. There exist $\varepsilon_0 >0$ and $l\in \N $ such that if
\begin{itemize}
\item $\varepsilon <\varepsilon_0$,  $S \in {\mathcal A}(M,\varepsilon,\alpha )$ and $T$ is an $m$-block of $C_S$;
\item $F:T\to L$ is a  $ C^{1}$  map with $\mathcal{N}(F)<1/2 $;
\item the image by $F$ of a connected component of $C^{c}_S$ is a  connected component of $C^{c}_S$;
\item there exists $J_i(m+p)\subset T$ such that $F( J_i(m+p))=J_j$;
\end{itemize}
 then $p\leq l$.
 \end{lemma}

\dem Let $\mu_i=\displaystyle\frac{|L|}{|I_i|}$, $1\leq i\leq k$. Let $\delta >0$ such that
$$\min \{\mu_i-\delta :\ \  1\leq i\leq k\} >1$$
and let $\beta = \min \{\mu_i-\delta :\ \ 1\leq i\leq k\}$.

Let $\varepsilon_0$ such that if $\varepsilon <\varepsilon_0$, $S\in {\mathcal A}(M,\varepsilon,\alpha )$ and $T$ is
a $m$-block of $C_S$, then  $|S_{\mid_{I_{i}}}'|>\beta $ and (see Remark \ref{rk1}) $\mathcal{N}(S^{m}_{\mid_{ T}})<1/2$.

Now, for an  $m$-block $T$ of $C_S$, and $F$ as in the hypotheses of the lemma, if $J_i(m)\subset T$, then

\begin{equation}\label{eq9}
\frac{|J_i(m+p)|}{|J_i(m)|}      \exp(-\frac{1}{2})  \stackrel{\mbox{\scriptsize
\begin{tabular}{c}
$\mathcal{N}(F)<\frac{1}{2} $
\end{tabular}}}\leq \frac{|F(J_i (m+p))|}{|F(J_i (m))|}  \stackrel{\mbox{\scriptsize
\begin{tabular}{c}
 $\mathcal{N}(F)<\frac{1}{2} $
\end{tabular}}}  \leq\exp(\frac{1}{2})\frac{|J_i(m+p)|}{|J_i(m)|}
\end{equation}

On the other hand,

\begin{equation}\label{eq10}
    \frac{|J_i(p)|}{|J_i|} = \frac{|S^{m}(J_i(m+p))|}{|S^{m}(J_i(m))|} \stackrel{\mbox{\scriptsize
\begin{tabular}{c}
 $\mathcal{N}(S^{m}_{\mid_{ T}})<\frac{1}{2} $
\end{tabular}}}   \geq  \exp(-\frac{1}{2})\frac{|J_i(m+p)|}{|J_i(m)|}.
  \end{equation}

Since $|S_{\mid_{I_{i}}}^{'}|>\beta $,

\begin{equation}\label{eq1111}
  |J_i(p)|\leq \frac{|J_i|}{\beta ^{p}}
\end{equation}

By inequalities (\ref{eq9}),  (\ref{eq10}) and  (\ref{eq1111}),

\begin{equation}\label{eq11}
 \frac{|F(J_i (m+p))|}{|F(J_i (m))|}
   \leq  \frac{\exp(1)}{\beta^{p}}.
\end{equation}

Since $F(J_i (m+p)))= J_j$,

\begin{equation}\label{eq11}
  |F(J_i (m))|\geq \frac{\beta^{p} |J_j|}{\exp(1)}.
\end{equation}

Consequently, because $|F(J_i (m))|\leq \max\{|J_r|\mid   1\leq r\leq k-1\}$ and $\beta >1$,
 there exists $l\in \mathbb{N}$ such that $0\leq p\leq l$.\hfill $\square$

\begin{corollary}\label{lema3}
Given $\alpha >0 $, $M >0 $ and $\eta$, with $0<\eta <1/2$, let $\varepsilon$ be as in
the  Lemma \ref{lema2}.
There exists a natural number $l $ such that if $S \in {\mathcal A}(M,\varepsilon,\alpha )$, and $C_S$ is $C^1$-minimal for some diffeomorphism $f_{_{S}}: S^1\to S^1$,
there exists $m_0=m_0(f_{_{S}})$ such that if $T$ is an $m$-block, with $m>m_0$ and
$J_i(m+p)\subset T$ with $ g_{_T}(J_i (m+p))= J_j$, then $0\leq p\leq l$.
\end{corollary}

\vspace{0.5cm}

 We denote  $\theta_p$ a $p$-uple $(i_1, \dots, i_p) \in \{1,\dots, k\}^p$ if $p>0$ and  $\theta_0 =0$ for $p=0$.

  Fix $M,\varepsilon,\alpha$, and take  $S \in {\mathcal A}(M,\varepsilon,\alpha )$ with
   $C_S$  minimal for a $C^1$ diffeomorphism $f:S^1\to S^1$. The  Corollary \ref{lema3}
   say that if $\varepsilon$ is small enough, the set $E_S$ of triples  $(i, j, \theta_p)$, $\theta_p =(i_1, \dots, i_p)$, for which there is an $m$-block $T= \phi^S(j_1, \dots, j_m)(L)  $ of $C_S$, with $m>m_0(f)$ satisfying
 \begin{equation}\label{triple2}
g_{_T}(J_i (m+p))= J_j,
\end{equation}
  with $J_i(m+p)= \phi^S(j_1, \dots, j_m,i_1, \dots, i_p)(J_i) \subset T$,  is a finite set.

   Remark that if $T$ is an $m$-block and  $g_{_T}(J_i (m))\neq  J_j$ for any $1\leq i,j\leq k-1$, then we can take the smaller $r\geq m$ such that $g_{_T}(J_i (r))= J_j$ for some $1\leq i, j\leq k-1$ and we can replace $T$ by the $r$-block contained in $T$ that contains this $J_i(r)$.  It follows that the equality (\ref{triple2}) is true for at least one  triple $  (i,j,\theta_0)=(i, j, 0)$.

 We call $B_S$  the subset of $E_S$ containing the triples
 $(i,j, \theta_p)$  such that the equality (\ref{triple2}) is true for infinite many blocks $T$ of $C_S$. We name $\mathcal{T}_{B_{S}}$ the set of such blocks and $\mathcal{T}_{(ij \theta_p)_S}$ the set of the blocks in $\mathcal{T}_{B_{S}}$ associated to $(i,j,\theta_p)$.

 Remark that $B_S\neq \emptyset$, and it is possible to cover $C_S$ only by blocks in $\mathcal{T}_{B_{S}}$. More than that, being  $E_S$  finite, by Corollary \ref{lema3}, there is an $m_0=m_0(f)$ such that for any $m>m_0$, every $m$-block $T$ is in $\mathcal{T}_{B_{S}}$.

In the following Remark we want to express the fact that if $ S \in {\mathcal A} (M, \varepsilon, \alpha) $ with $\varepsilon$ sufficiently small, $S$ is $C ^ 1$ close to  the affine map  $ R\in  {\mathcal A}(M,\varepsilon,\alpha ) $ that satisfies $R'(x)S'(x) >0$. Also, fixed $ l\in \mathbb{N} $,  the measure of $J^{R}_j(i_1,....i_r)$ is close to the measure of $J^{S}_j(i_1,...,i_r)$  for $r\leq l$ and for every $ 1 \leq j \leq k- 1$.

\begin{rk} \label{c1closer}
Let $l \in \mathbb{N}$ and $M>0$, $\alpha>0$. Given $c>0$, there exists $\varepsilon>0$ such that if $S\in {\mathcal A}(M,\varepsilon,\alpha )$, then:

 $$|J^{R}_j(i_1,....i_r)| \exp(-c)  \leq |J^{S}_j(i_1,....i_r)|\leq \exp(c) |J^{R}_j(i_1,....i_r)|    $$
  for $r \leq l$ and $ 1 \leq j \leq k- 1$. It follows because $S$ is $C^1$ close to $R$ if $\varepsilon $ is close to $0$.

It follows
\begin{equation} \label{eqc1closer}
\exp(-2c)\frac{|J^{R}_j(i_1,....i_r)|}{| J^R_i(j_1,...,j_s)|}\leq  \frac{|J^{S}_j (i_1,....i_r)|}{|J^{S}_i (j_1,...,j_s)|} \leq \exp(2c)\frac{|J^{R}_j(i_1,....i_r)|}{| J^R_i(j_1,...,j_s)|},
\end{equation}
for all $1\leq i, j\leq k-1$ and $r, s\leq l$.

\end{rk}

 \begin{lemma}\label{lema4}
Let $M>0$ and $\alpha >0$. Suppose that there exist a sequence of positive numbers $\{\varepsilon_n\}$, $\varepsilon_n\to 0$, a sequence $\{S_n\}$,  $S_n\in  {\mathcal A}(M,\varepsilon_n,\alpha )$, with each $C_{S_{n}}$ $C^1$ minimal for some $f_{S_{n}}$ and $S'_r(x)S'_s(x) >0$ for all $x\in \cup I_i$ and $r,s\in \mathbb{N}$. Also suppose there exists a triple $(i,j,\theta_p)$ such that $(i,j,\theta_p) \in B_{S_{n}}$ for all $n \in \mathbb{N}$ and let $R$ be  the affine map in $  {\mathcal A}(M,\varepsilon,\alpha ) $ with $R'(x)S'_n(x) >0$ for all $x\in \cup I_i$.
Under such conditions, there exists an affine map  $h_{ij\theta_p}:  L \to L$, such that
 $h_{ij\theta_p}(C_R)\subset  C_R$, $h_{ij\theta_p}(C^c_R)\subset  C^c_R$ and $h_{ij\theta_p}(J_i(p))=J_j$ where $\theta_p =(i_1, \dots, i_p)$ and $J_i(p)=\phi^R(i_1, \dots, i_p)(J_i)$.
\end{lemma}

 \dem To make the argument clear, first we prove the Lemma for a triple $(i,j,\theta_0)=(i,j,0)$.

{\bf{Claim}}:  Given $J^{R}_l(r)\subset L$, with $1\leq l \leq k-1$ and $r \in \N$ there exists a connected component $J$ of $ C_R^{c}$ such that  $|J|=\frac{|J_j|}{| J_i|} |J^{R}_l(r)|$.

By Lemma \ref{lema2} and Remark \ref{rk1}, for all $\eta >0$, there exists
  $n_0=n_0 (\eta)$ such that for $n\geq n_0$, if $S_n\in {\mathcal A}(M,\varepsilon_n,\alpha )$
 with  $C_{S_{n}}$  minimal for a $C^1$-diffeomorphism $f_{S_{n}}$, and $T^{S_n}$ is  an $m$-block for $S_n$ with $m>m_0(f_{S_{n}})$, then
 \begin{equation} \label{eqeta}
 \mathcal{N}({g_{_{T^{S_n}}}}) < \eta \ \mathrm{and} \
\mathcal{N}({S^{m}_{n}}_{\mid_{T^{S_n}}}) < \eta.
\end{equation}

 Let $\eta>0$,  $n> n_0(\eta)$ and $m>m_0(f_{S_{n}})$, as above, and let $T^{S_n}=\phi^{S_{n}}(j_1,...,j_m)$ be an $m$-block in $\mathcal{T}_{(ij \theta_0)_{S_n}}$. Given $J^{R}_l(r)=\phi^R(i_1, \dots, i_r)(J_l)$, for each $n\in \mathbb{N}$ consider $J^{S_{n}}_l(r)=\phi^{S_n}(i_1, \dots, i_r)(J_l)$. Let $J^{S_{n}}_l(m+r)$ be $\phi^{S_n}(j_1, \dots, j_m)(J^{S_{n}}_l(r))$ contained
  in $T^{S_n}$. Then,  by inequalities  (\ref{eqeta}) and reasoning as in the inequalities (\ref{eq9}) and (\ref{eq10}),

\begin{equation}\label{eq323}
  \exp(-2\eta)\frac{|J^{S_{n}}_l(r)|}{| J_i|}\leq  \frac{|g_{_{T^{S_n}}}(J^{S_{n}}_l (m+r))|}{|g_{_{ T^{S_n} }}(J^{S_{n}}_i (m))|} \leq \exp(2\eta)\frac{|J^{S_{n}}_l(r)|}{| J_i|}.
  \end{equation}

Since $g_{_{  T^{S_n}}}(J^{S_{n}}_i (m))=  J_j$, it follows

\begin{equation}\label{eq324}
 \exp(-2\eta )\frac{|J^{S_{n}}_l(r)||J_j|}{| J_i|}\leq  |g_{_{T^{S_n}  }}(J^{S_{n}}_l (m+r))| \leq \exp( 2\eta )\frac{|J^{S_{n}}_l(r)||J_j|}{| J_i|}.
\end{equation}

If $\varepsilon_n$ is small enough, $S_n$ is close to $R$ in the $C^1$ topology (see Remark \ref{c1closer}) and then $|J^{S_{n}}_l(r)|$  can be made to be as close to  $|J^{R}_l(r)|$ as desired by taking $n$ big enough.  Therefore,
if $\eta$ is small enough and $n$ is great enough (it is means $\varepsilon_n$ sufficiently small),
there is  only a finite number of possibilities for $|g_{_{ T^{S_n} }}(J^{S_{n}}_l (m+r))|$ for each $n >n_0$ big enough and any $m> m_0(f_{S_{n}})$.

 Now we will prove that there exist an unique $\theta' _{r'}=(i'_1,...,i'_{r'})$ and an unique $l^{'}$ such that
$$g_{_{  T^{S_n}}}(J^{S_{n}}_l (j_1,...,j_m, i_1,\dots, i_r))   = \phi^{S_{n}}(i'_1,...,i'_{r^{'}}) (J_{l^{'}} )$$
 for $n\geq n_1$, for a suitable $n_1$ and $m> m_0(f_{S_{n}})$.

Let $\eta >0$ and  $n>n_0(\eta)$ and $m_0$ as in the inequalities  (\ref{eqeta}). Let  $  T^{S_n} = \phi^{S_n}(j_1,...,j_m)$ be an $m$-block of $C_{S_n}$  with $m > m_0$, contained in $\mathcal{T}_{(ij \theta_0)_{S_n}}$ and $J^{S_{n}}_l (j_1,...,j_m, i_1,\dots, i_r)\subset T^{S_n}$. Denote $g_{_{ T^{S_n} }}(J^{S_{n}}_l (j_1,...,j_m, i_1,\dots, i_r))$ by $J^{S_n}$. By inequality (\ref{eq324}), $|J^{S_n}|$ is as close to $|J^R_l(r)|\frac{|J_j|}{|J_i|}$ as we want if $n$ is big enough.

\begin{figure}
\psfrag{gt}{$g_{_{T^{S_{n}}}}$}\psfrag{r}{$S_n^{m}$}
\psfrag{tm}{$T^{S_{n}}$}\psfrag{l}{$L$}\psfrag{jj}{$J_j$}
\psfrag{j}{$J^{S_{n}}=g_{_{T^{S_{n}}}}(J^{S_{n}}_l(m+r))$ }
\psfrag{jl}{$J^{S_{n}}_l(r)$}
\psfrag{jp}{$J^{S_{n}}$}
\psfrag{c2}{$c_2(n)$}\psfrag{c1}{$c_1(n)$}\psfrag{cm}{$c_1(m,n)$}
\psfrag{jim}{$J^{S_{n}}_i(m)$}
\psfrag{jlm}{$J^{S_{n}}_l(m+r)$}
\psfrag{ji}{$J_i$}
\psfrag{jlr}{$J^{S_{n}}_l(r)$}

\psfrag{jii}{$J_i$}\psfrag{jlrr}{$J^{R}_l(r)$}\psfrag{c11}{$c_1$}
\psfrag{rr}{$R$}

\begin{center}
\caption{\label{figura4}}
{\includegraphics[scale=0.25]{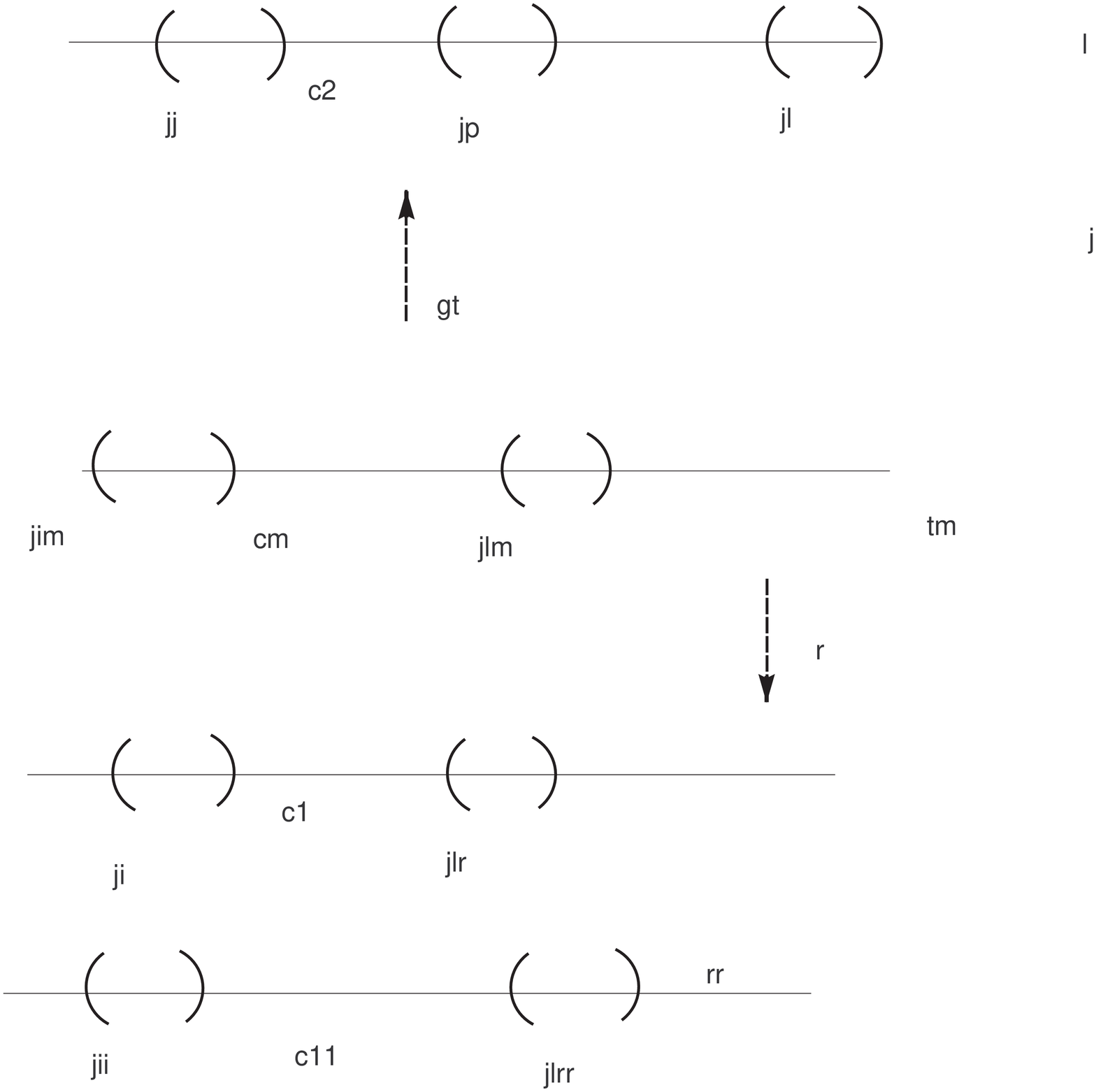}}
\end{center}
\end{figure}

Let $c_1$, $c_1(n)$, $c_1(m,n)$ and $c_2(n)$ the lengths of the intervals indicated in the Figure \ref{figura4}. Then, because $S_n^{m}(J^{S_{n}}_i(m))=J_i$,

\begin{equation}\label{eq555}
  \frac{|J_i|}{c_1(n)} \stackrel{\mbox{\scriptsize
\begin{tabular}{c}
$ x,y\in T^{S_n}  $
\end{tabular}}}   =   \frac{|J^{S_{n}}_i(m)|   (S_n^{m})^{'}(x)}{c_1(m,n) (S_n^{m})^{'}(y) } \stackrel{\mbox{\scriptsize
\begin{tabular}{c}
$ \mathcal{N}({S^{m}_{n}}_{\mid_{T^{S_n} }}) < \eta$
\end{tabular}}} \leq \exp( \eta)  \frac{|J^{S_{n}}_i(m)| }{c_1(m,n) }.  \end{equation}
Therefore,

\begin{equation}\label{eq888}
  \frac{|J_i|}{c_1(n)} \exp(-\eta) \leq   \frac{|J^{S_{n}}_i(m)| }{c_1(m,n) }\leq  \exp(\eta)  \frac{|J_i|}{c_1(n)}.
  \end{equation}

On the other hand, because $g_{_{ T^{S_n} }}(J^{S_{n}}_i(m))=J_j$,
\begin{equation}\label{eq666}
   \frac{c_2(n)}{|J_j|} \stackrel{\mbox{\scriptsize
\begin{tabular}{c}
$ x\in  T^{S_n} $
\end{tabular}}}  = \frac{c_1(m,n). |g_{_{ T^{S_n} }}'(x)|}{|g_{_{ T^{S_n} }}(J^{S_{n}}_i(m))|}\stackrel{\mbox{\scriptsize
\begin{tabular}{c}
$ y\in  T^{S_n} $
\end{tabular}}}  = \frac{c_1(m,n).| g_{_{ T^{S_n}  }}'(x)|}{|J^{S_{n}}_i(m)|| g_{_{ T^{S_n} }}'(y)|} \leq \exp(\eta) \frac{c_1(m,n)}{|J^{S_{n}}_i(m)|}.
\end{equation}
Therefore,

\begin{equation}\label{eq999}
\exp(-\eta) \frac{c_1(m,n)}{|J^{S_{n}}_i(m)|} \leq \frac{c_2(n)}{|J_j|} \leq \exp(\eta) \frac{c_1(m,n)}{|J^{S_{n}}_i(m)|}.
\end{equation}

From the equalities (\ref{eq888}) and (\ref{eq999}), we conclude

 \begin{equation}\label{eq777}
  \exp(-2\eta)  \frac{|J_j|}{|J_i|}c_1(n)  \leq c_2(n) \leq   \exp(2\eta) \frac{|J_j|}{|J_i|}c_1(n).
\end{equation}

Therefore, since $l$ and $r$ are fixed, so $c_1(n)\to c_1 $ when $n\to \infty$ and then, from inequalities (\ref{eq777}),  $c_2(n) \to \frac{|J_j|}{|J_i|}c_1$. So, not only there is a finite number of possibilities for  the length of $g_{_{ T^{S_n} }}(J^{S_{n}}_l (j_1,...,j_m, i_1,\dots, i_r))$, as there is a finite number of possibilities  for the connected component of $C^{c}_{S_n}$ which is image of $J^{S_{n}}_l (j_1,...,j_m, i_1,\dots, i_r)$ by $g_{_{T^{S_n}  }}$,  for all $n$ and $m$ big enough.

If $g_{_{T^{S_n}}}$ is increasing or decreasing for all $m$-blocks $T^{S_n}$ for $n$ great enough and for all  $m> m_0(f_{_{S^n}})$, then  we conclude there exist an unique $l^{'}$ and an unique $\theta _{r'}=(i'_1,...,i'_{r'})$ such that
$g_{_{ T^{S_n} }}(J^{S_{n}}_l (j_1,...,j_m, i_1,\dots, i_r))    = \phi^{S_{n}}(i'_1,...,i'_{r'}) (J_{l'} )$ for $n\geq n_1$, for a suitable $n_1$ and $m> m_0(f_{_{S_n}})$.

Because  $\phi^{S_{n}}(i'_1,...,i'_{r'}) (J_{l'} )$ goes to $\phi^{R}(i'_1,...,i'_{r'}) (J_{l'} )$ in the Hausdorff distance and  $|J^{S_n}|\to  |J^R_l(r)|\frac{|J_j|}{|J_i|}$ then $|\phi^{R}(i'_1,...,i'_{r'}) (J_{l'} )  | =|J^R_l(r)|\frac{|J_j|}{|J_i|}.$

Let $h_{ij\theta_0}:L \to L$ the unique affine map such that $h_{ij\theta_0}(J^R_l(r)) =  \phi^{R}(i'_1,...,i'_{r'}) (J_{l'} )$. Note that $h_{ij\theta_0}$ verifies the following  properties:

 \begin{enumerate}
 \item $h_{ij\theta_0}(J_i)= J_j$
 \item  $h_{ij\theta_0}(C_R)\subset C_R$ and $h_{ij\theta_0}(C_R^c)\subset C_R^c$.
 \item  $|h_{ij\theta_0}(J^R_l(r))|=\frac{|J_j|}{|J_i|}|J^R_l(r)|$ for any $1\leq l\leq k-1$ and any $r\in \mathbb{N}$.
 \end{enumerate}

If for all natural numbers $M$ and $N$ there exist $n >N$ and $m_1, m_2> M$,  $T_1^{S_n}$ is an $m_1$-block  and $T_2^{S_n}$ is an $m_2$-block of $C_{S_n}$ such that $g_{_{T_1^{S_n}}} $ is an increasing map  and $g_{_{T_2^{S_n}}} $ is an decreasing map, with an analogous reasoning, we determine exactly two different affine  maps $h_{ij\theta_0}$ for the triple $(i,j, \theta_p)$.

The proof for the case $(i,j, \theta_p)$ with $p>0$ is very similar, replacing $ J^{S_{n}}_i(m)$ by  $ J^{S_{n}}_i(m+p)$ in the equations used in this proof.
\hfill $\square$

\begin{figure}[h]
\psfrag{0}{$0$}\psfrag{1}{$\frac{1}{4}$}
\psfrag{2}{$\frac{3}{8}$}
\psfrag{3}{$\frac{1}{2}$}
\psfrag{4}{$\frac{3}{4}$}
\psfrag{5}{$1$}
\psfrag{9}{$0$}
\psfrag{10}{$\frac{1}{5}$}
\psfrag{11}{$\frac{2}{5}$}
\psfrag{12}{$\frac{3}{5}$}
\psfrag{13}{$\frac{4}{5}$}
\psfrag{14}{$1$}
\psfrag{i2}{$\frac{1}{3}$}\psfrag{i3}{$\frac{1}{2}$}\psfrag{l}{$R$}
\psfrag{f}{$C_R$\mbox{ is the usual ternary Cantor} }
\begin{center}
\caption{\label{figura1}}
\subfigure[]{\includegraphics[scale=0.1535]{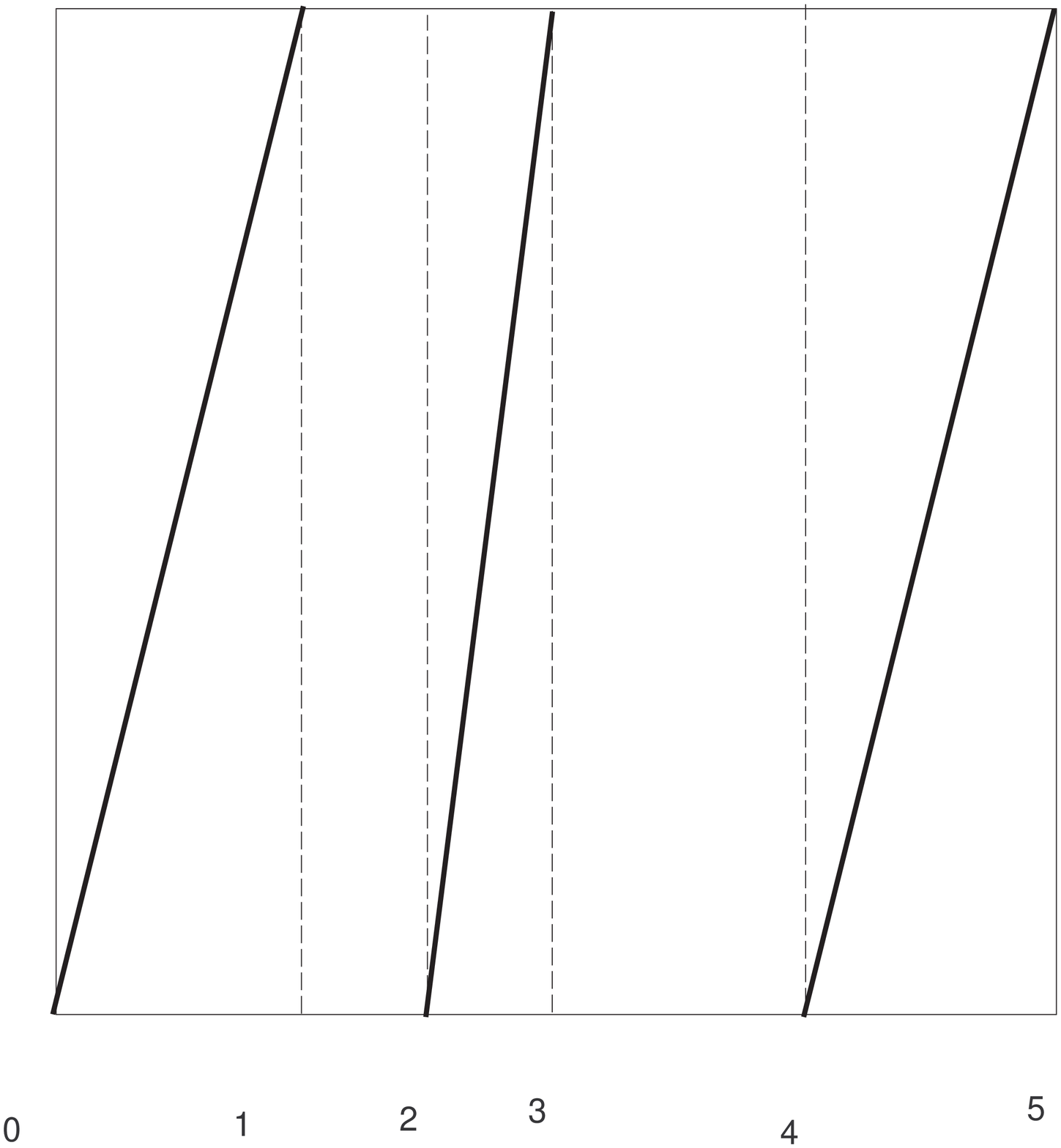}}
\subfigure[]{\includegraphics[scale=0.24]{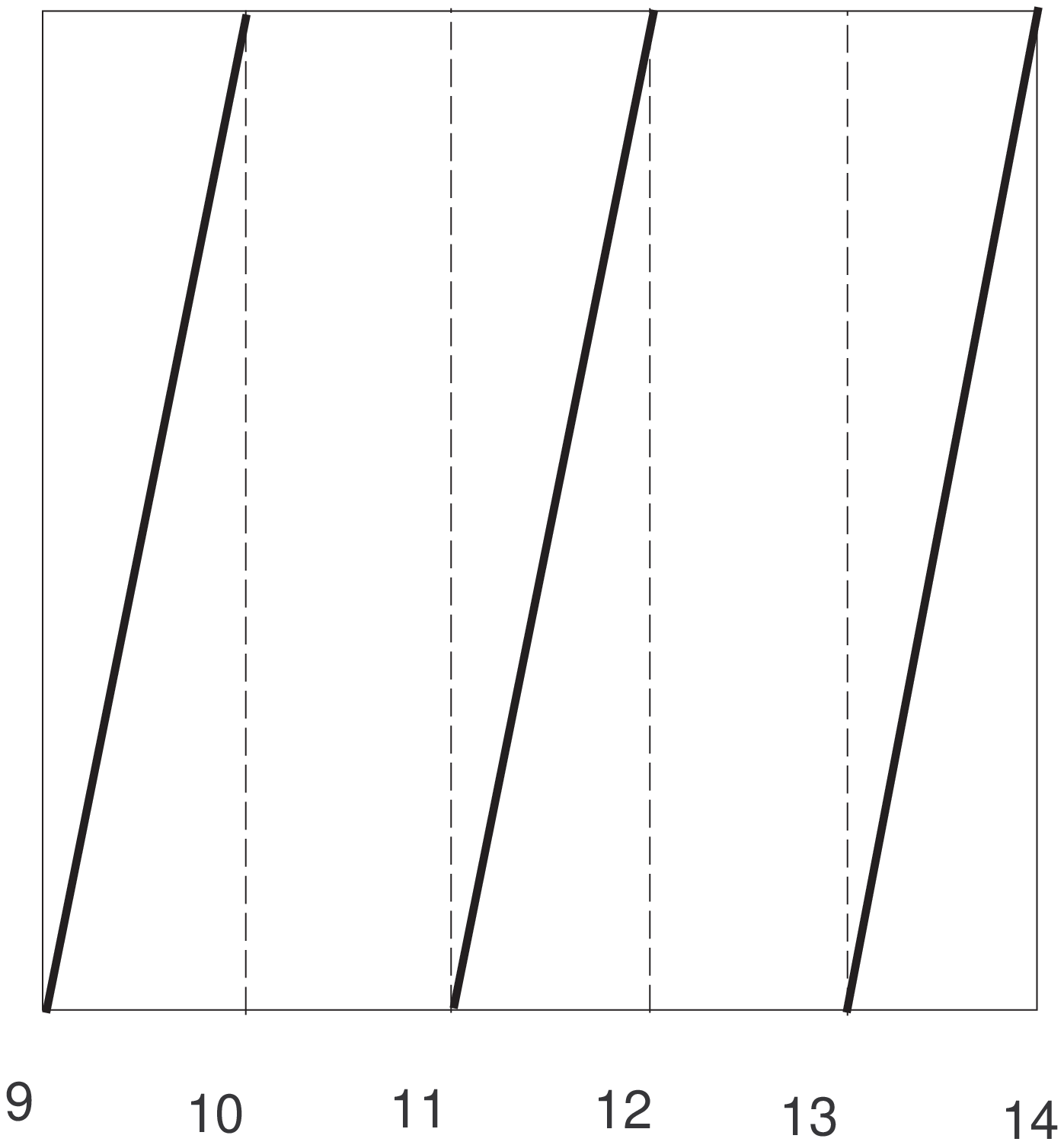}}
\end{center}
\end{figure}
In the Figure \ref{figura1}(a), we exhibit an example in which the set $A_{ij\theta_p}$ of the previous Lemma can be the interval $[0,1/2]$ and in the Figure \ref{figura1}(b) we exhibit an example in which the only possibility for the set $A_{ij\theta_p}$ is the entire interval $[0,1]$.

\begin{lemma}\label{lema5} Let $l \in \mathbb{N}$ and $M>0$, $\alpha> 0$. There exist $\eta $ with $0< \eta <1/2$ and $\varepsilon >0$ such that if
\begin{itemize}
\item $S, Z \in {\mathcal A}(M,\varepsilon,\alpha )$;
\item $T^S$ is an $m$-block of $C_S$, $T^S=\phi (i_1,...,i_m)(L) $ and $T^Z$ is an $m'$-block of $C_Z$ with $T^Z=\phi (i'_1,...,i'_{m'})(L) $;
\item $F:T^S\to L$ and $G:T^Z\to L$ are $C^1$ maps such that $\mathcal{N}(F) <\eta$, $\mathcal{N}(G) <\eta$, the image by $F$ (by $G$) of a connected component of $C^{c}_S$ (of $C^{c}_Z$) is a  connected component of $C^{c}_S$ (of $C^{c}_Z$) and $sig(F')=sig(G')$;
\item $F(\phi^S(i_1,\dots i_{m}, j_1, \dots, j_p)(J_i)) =J_j = G(\phi^Z(i'_1,\dots i'_{m'}, j_1, \dots, j_p)(J_i))$ for some $\theta_p=   ( j_1, \dots, j_p)$, $p\leq l$ and for some $i$ and $j$ with  $1\leq i, j\leq k-1$,
\end{itemize}
then for all $\theta_r=(q_1,...,q_r)$ with $0\leq r\leq 2l+2 $, $F(J^S_t(i_1,...,i_m,q_1,...,q_r))= \phi^S(q'_1,\dots, q'_{r'})(J_{t'})$ and $G(J^Z_t(i'_1,...,i'_{m'},q_1,...,q_r))=\phi^Z(q'_1,\dots, q'_{r'})(J_{t'})$ for some $r' \in \mathbb{N}$  and some $t'$ with $1\leq t'\leq k-1$.
\end{lemma}

\dem  Let $M>0$ and $\alpha> 0$. Given $\eta>0$, by Remark \ref{rk1} there exists $\varepsilon >0$  such that if $S \in {\mathcal A}(M,\varepsilon,\alpha )$ and  $T^S$  is an $m$-block of $C_S$, then
\begin{equation}\label{1010}
\mathcal{N}({S^{m}}_{\mid_{T^S}}) < \eta.
\end{equation}

Let $F$ be as in the hypotheses of the Lemma, with $\mathcal{N}(F) <\eta$. Let $T^S=\phi^S(i_1,\dots i_{m})(L)$. Consider $J^S_i(i_1,...,i_m,j_1,...,j_p)$ with $p\leq l$ and $J^S_t(i_1,...,i_m,q_1,...,q_r)$ with $r\leq 2l+2$. Under these conditions,

\begin{equation}\label{eq3333}
\exp(-\eta)\frac{|J^S_i(i_1,...,i_m,j_1,...,j_p)|}{|J^S_t(i_1,...,i_m,q_1,...,q_r)|}\leq \frac{|F(J^S_i(i_1,...,i_m,j_1,...,j_p))|}{|F(J^S_t(i_1,...,i_m,q_1,...,q_r))|}\leq \exp(\eta)\frac{|J^S_i(i_1,...,i_m,j_1,...,j_p)|}{|J^S_t(i_1,...,i_m,q_1,...,q_r)|}\nonumber
\end{equation}
so, by the inequality (\ref{1010}), because $S^{m}(J^S_i(i_1,...,i_m,j_1,...,j_p))=  J^S_i(j_1,...,j_p)$ then

\begin{equation}
\exp(-2\eta)\frac{|J^S_i(j_1,...,j_p)|}{|J^S_t(q_1,...,q_r)|}\leq \frac{|F(J^S_i(i_1,...,i_m,j_1,...,j_p))|}{|F(J^S_t(i_1,...,i_m,q_1,...,q_r))|}\leq \exp(2\eta)\frac{|J^S_i(j_1,...,j_p)|}{|J^S_t(q_1,...,q_r)|}\nonumber
\end{equation}
Therefore, by Remark \ref{c1closer},
\begin{equation}
\exp(-3\eta)\frac{|J^R_i(j_1,...,j_p)|}{|J^R_t(q_1,...,q_r)|}\leq \frac{|F(J^S_i(i_1,...,i_m,j_1,...,j_p))|}{|F(J^S_t(i_1,...,i_m,q_1,...,q_r))|}\leq \exp(3\eta)\frac{|J^R_i(j_1,...,j_p)|}{|J^R_t(q_1,...,q_r)|}\nonumber
\end{equation}
and, since $F(J^S_i(i_1,...,i_m,j_1,...,j_p)) =J_j$, it follows

\begin{equation}\label{eqfg}
\exp(-3\eta)\frac{|J_j||J^R_t(q_1,...,q_r)|}{|J^R_i(j_1,...,j_p)|}\leq {|F(J^S_t(i_1,...,i_m,q_1,...,q_r))|}\leq \exp(3\eta)\frac{|J_j||J^R_t(q_1,...,q_r)|}{|J^R_i(j_1,...,j_p)|}.
\end{equation}
Analogously for $Z$ and $G$ it follows
\begin{equation}\label{eqfg1}
\exp(-3\eta)\frac{|J_j||J^R_t(q_1,...,q_r)|}{|J^R_i(j_1,...,j_p)|}\leq {|G(J^Z_t(i^{'}_1,...,i^{'}_{m^{'}},q_1,...,q_r))|}\leq \exp(3\eta)\frac{|J_j||J^R_t(q_1,...,q_r)|}{|J^R_i(j_1,...,j_p)|}.
\end{equation}

Now, with an  analogous  argument used in the proof of Lemma \ref{lema4}, we conclude that if $\eta >0$ is token small enough there is $\varepsilon >0$ such that remains only one choice for the image by $F$ for $J^S_t(i_1,...,i_m,q_1,...,q_r)$, with $p\leq l$ and $r\leq 2l+2$. It is means, $F(J^S_t(i_1,...,i_m,q_1,...,q_r))= J^S_{t'}(q'_1,...,q'_{r'})$  for some $t'$ and $(q'_1,...,q'_{r'})$ that do not depend on $S$. This conclude the proof. \hfill$\square$

\begin{lemma}\label{lema6}
Under the same hypotheses of  Lemma \ref{lema4}, there exists $n_0\in \mathbb{N}$ and, for all $n>n_0$, there exists  $m_0 = m_0(f_{S_n})$ such that if
\begin{itemize}
\item $n>n_0$;
\item  $m> m_0$;
\item  $T^{S_n}$ is an $m$-block of $C_{S_n}$, $T^{S_n}=\phi^{S_{n}}(i_1,\dots i_{m})(L)$, $T^{S_n}\in \mathcal{T}_{(ij \theta_p)_{S_n}}$;
\item  $g'_{_{T^{S_n}}}$ has the same signal of $g'_{_{I^{S_v}}}$ for infinitely many $I^{S_v}\in \mathcal{T}_{(ij \theta_p)_{S_v}}$ for infinitely many natural numbers $v$
\end{itemize}
 then, for each $(q_1,...,q_r)$ and $t$, $1\leq t \leq k-1 $, there exist $(q'_1,...,q'_{r'})$ and $t'$, $1\leq t'\leq k-1 $ such that $g_{_{T^{S_n}}}(J^{S_{n}}_t( i_1,\dots i_{m}, q_1,...,q_r))=J^{S_n}_{t'}(q'_1,...,q'_{r'})$, where $t'$ and $(q'_1,...,q'_{r'})$ do not depend on $n$ neither on  $T$, but can depend on the signal of $g'_{_{T^{S_n}}}$.
\end{lemma}

\dem  By Lemma \ref{lema2}, there exists $n_1\in \mathbb{N}$ and,  for each $n>n_1$, if  $C_{S_n}$ is minimal for a $C^ 1$ diffeomorphism  $f_{_{S_n}}: S^ 1\to S^ 1$, there exists $m_0(f_{_{S_n}}) \in \mathbb{N}$ such that if  $T^{S_n}$ is a $m$-block of $C_{S_n}$,
with $m> m_0(f_{_{S_n}})$,  then
$$\mathcal{N}({g_{_{T^{S_n}}}}) < 1/2.$$

Let  $l$ given by Corollary \ref{lema3}. For such $l$, let $\eta >0 $ as assured in Lemma \ref{lema5}. Let $n_0 >n_1$ such that if $n> n_0$

$$\mathcal{N}( S^{m}_n|_{T^{S_n}})<\eta/2 $$
for all $m$-block $T^{S_n}$ of $C_{S_n}$ (assured by Remark \ref{rk1}) and, also,

$$\mathcal{N}( g_{_{T^{S_n}}})<\eta/2 $$
 for all $m$-block $T^{S_n}$ of $C_{S_n}$ with $m>m_0(f_{_{S_n}})$ as in Lemma \ref{lema2}.

Let $a, b$ be natural numbers with $a, b\geq n_0$ and
$c\geq m_0( f_{_{S_a}})$ and $d\geq m_0( f_{_{S_b}})$.

Consider $T^{S_a} = \phi^{S_{a}} (i_1,...,i_{c})(L)$ in $\mathcal{T}_{(ij \theta_p)_{S_{a}}} $  and $T^{S_b} = \phi^{S_{b}} (i^{'}_1,...,i^{'}_{d})(L) $,
in  $\mathcal{T}_{(ij \theta_p)_{S_{b}}} $.

 So there exist $J_i^{S_{a}}(i_1,...,i_{c},j_1,...,j_p)$ and  $J_i^{S_{b}}(i^{'}_1,...,i^{'}_{d},j_1,...,j_p)$ such that
$$g_{_{T^{S_a}}}(J_i^{S_{a}}(i_1,...,i_{c},j_1,...,j_p))= g_{_{T^{S_b}}}(J_i^{S_{b}}(i^{'}_1,...,i^{'}_{d},j_1,...,j_p)) =J_j.  $$

Then by  Lemma \ref{lema5} taking $F=g_{_{T^{S_a}}}$ and $G= g_{_{T^{S_b}}}$ , for each $t$, $1\leq t\leq k-1$ and each $(q_1,...,q_r)$ with $r \leq 2l+2$, there exist $t'$ and $(q'_1,...,q'_{r'})$ such

$$g_{_{T^{S_a}}}(J_t^{S_{a}}(i_1,...,i_{c},q_1,...,q_r))=J^{S_a}_{t'}(q'_1,...,q'_{r'}) \mbox{ and }$$
$$ g_{_{T^{S_b}}}(J_t^{S_{b}}(i^{'}_1,...,i^{'}_{d},q_1,...,q_r))=J^{S_{b}}_{t'}(q'_1,...,q'_{r'}).$$

Until here, we have proved the Lemma for $r\leq2l+2$. Now, we will prove the assertion for any $r \in \N$.

Let $ T^{S_a}(j_1) \subset T^{S_a}$ and $ T^{S_b}(j_1) \subset T^{S_b}$ . Consider $k_1,k_2\in \N$ such that

$$ k_1 = \max \{ k\in \mathbb{N}:\  S_{a}^{k-1}\circ g{_{_{T^{S_{a}}}}}_{\mid _{  T^{S_a}(j_1)}} \subset \cup I_i \},  $$

$$ k_2 = \max \{ k\in \mathbb{N}:\  S_{b}^{k-1}\circ g{_{_{T^{S_b}}}}_{\mid _{ T^{S_b}(j_1)}} \subset \cup I_i \}.  $$

 {\bf{Claim:}} $k_1=k_2$.

By definition of $k_1$ there exists  $J^{S_{a}}_s(i_1,...,i_c, j_1,u_1,...,u_r )$ such that
$$S_{a}^{k_{1}}\circ g_{_{T^{S_{a}}}}(J^{S_{a}}_s(i_1,...,i_c, j_1,u_1,...,u_r ) )=J_{t},$$
 for some $t$, $1\leq t\leq k-1.$

Since $\mathcal{N}(S_{a}^{k_1}) <\eta/2$ and $\mathcal{N}(g_{_{T^{S_a}}}) <\eta/2$, by {\cite[Lemma 1]{N}}, $$\mathcal{N}(  S_{a}^{k_1} \circ g_{_{T^{S_a}}}) <\eta.$$
 So, taking $F= S_{a}^{k_1} \circ g_{_{T^{S_a}}}$  by Lemma \ref{lema3}, $r\leq l$. Therefore, taking $F=g_{_{T^{S_a}}}$ and $G=g_{_{T^{S_b}}}$, by Lemma \ref{lema5}, there exist $s'$ and  $ ( u'_1,...,u'_{r'})$ such that
  $$ g_{_{T^{S_a}}} (   J^{S_{a}}_s(i_1,...,i_c, j_1,u_1,...,u_r )) = J_{s'}^{S_{a}}( u'_1,...,u'_{r'})$$
   and  $ g_{_{T^{S_b}}} (  J^{S_{b}}_s(i^{'}_1,...,i^{'}_d, j_1,u_1,...,u_r )  ) = J_{s'}^{S_{b}}( u'_1,...,u'_{r'}   ) $.
Then
 $$ S_{a}^{k_1}\circ g_{_{T^{S_a}}} (   J^{S_{a}}_s(i_1,...,i_c, j_1,u_1,...,u_r )) =  J_t =S_{b}^{k_1}\circ g_{_{T^{S_b}}}( J^{S_{b}}_s(i^{'}_1,...,i^{'}_d, j_1,u_1,...,u_r )  ). $$
 This imply that $k_2\leq k_1$. Analogously $k_1\leq k_2$, so $k_1=k_2$.

Also, if $ {S_{a}^{k_1}}_{\mid g_{_{T^{S_a}}} (T^{S_a}(j_1))}^{-1} =\phi^{S_a}(\tau_1, \dots, \tau_{k_1})$, then

\begin{equation}\label{eqk1k2}
 {S_{b}^{k_1}}_{\mid g_{_{T^{S_b}}} (T^{S_b}(j_1))}^{-1} =\phi^{S_b}(\tau_1, \dots, \tau_{k_1}).
 \end{equation}

Note that $ S_{a}^{k_1}\circ g_{_{T^{S_a}}}$ and $  S_{b}^{k_1} \circ g_{_{T^{S_b}}}$ are under the hypotheses of Lemma \ref{lema5}.
Then for each $s$, $1\leq s\leq k-1$, $j_1$ and $(q_1,...,q_r)$ with $r\leq 2l+2 $, there exist $t'$ and $(q'_{1},...,q'_{r'})$ such that
$$S_{a}^{k_1}\circ g_{_{T^{S_a}}} (   J^{S_{a}}_s(i_1,...,i_c, j_1,q_1,...,q_r ))= J_{t'}^{S_a}( q'_1,...,q'_{r'}  )$$
 and $ S_{b}^{k_1} \circ g_{_{T^{S_b}}}( J^{S_{b}}_s(i^{'}_1,...,i^{'}_{d}, j_1,q_1,...,q_r )  )  =J_{t'}^{S_{b}}( q'_1,...,q'_{r'}  )  $.

 Therefore, by (\ref{eqk1k2}), there exist $s''$ and $(w_1,...w_{r''})$ such that
 $$ g_{_{T^{S_a}}} (   J^{S_{a}}_s(i_1,...,i_c, j_1,u_1,...,u_r ))= J_{s''}^{S_{a}} ( w_1,...w_{r''}  )$$
  and
 $ g_{_{T^{S_b}}}( J^{S_{b}}_s(i^{'}_1,...,i^{'}_{d}, j_1,u_1,...,u_r )  ) = J_{s''}^{S_{b}} ( w_1,...w_{r''}  ) $
  for $r\leq 2l+2$.

Proceeding inductively, the Lemma follows. \hfill $\square$

\vspace{0.5cm}

\noindent {\bf{Proof of Theorem \ref{teorema}: }}  Let $M>0$ and $\alpha >0$. If there exists $\varepsilon >0$ such that  for all $S\in {\mathcal A}(M,\varepsilon,\alpha )$,   $C_S$ is not minimal for any $C^1$ diffeomorphism of $S^1$, then the Theorem holds.

So suppose that  for all $\varepsilon >0$ there exists $S \in {\mathcal A}(M,\varepsilon,\alpha )$ such that $C_S$ is $C^1$-minimal. Then, by Lemma \ref{lema4} and Lemma \ref{lema6}, we can find $\varepsilon>0$ and $ S \in {\mathcal A}(M,\varepsilon,\alpha )$ such that $C_S$ is minimal
 for a $C^1$-diffeomorphism $f_{_{S}}: S^1\to S^1$ and it is possible  to cover $C_S$ by  $T_1, \dots, T_n$, $T_l = \phi^S(j_1,\dots, j_{m_l})(L)$ for all $1\leq l\leq n$, $T_l\cap T_r=\emptyset$  if $l\neq r$, and
 if $R\in {\mathcal A}(M,\varepsilon,\alpha )$ is the affine map such that $S'(x)R'(x) >0$ for all $x\in \cup I_i$, for each $1\leq l\leq n$, there
exists an affine map $h_l:L\to L$, such that $h_l(C_R^c)\subset C_R^c$, $h_l(C_R)\subset C_R$ and ,  if  $h_{l}  ( J^R_j    ( q_1,...,q_r)   )= J^R_{t}(i_1,...,i_{s} )$, then
 $$ g_{_{T_l}}(  J_j ^{S}   ( j_1,...,j_{m_l},q_1,...,q_r)) =J^S_{t}(i_1,...,i_{s} ).$$

Since $g_{_{T_l}} = S^{u_{T_l}}\circ {f_{_{S}}}_{\mid_{T_l}}$,  there exist $(t_1,..., t_{u_{T_l}})$ such that $  {f_{_{S}}}_{\mid_{T_l}}=\phi^S(t_1, \dots, t_{u_{T_l}}) \circ g_{_{T_l}}$, for all $1\leq l\leq n$.

Now we will define a diffeomorphism $F: S^1 \to S^1$ of class $C^{\infty}$. Let $\psi$ be, as in the Lemma \ref{lema1}.

First define $F_{\mid_{\psi (T_l) \cap  C_R^{c}  }}$ such that
 $$ F(J^R_l(j_1,...,j_{m_l}, q_1,...,q_r)) =  \phi^R_{t_1}\circ...\circ \phi^R_{t_{u_{T_l}}} \circ h_{l} (  J_l^R (q_1,...,q_r) )    .$$
 for all $1\leq l\leq n$. Now,  extend this affine map $F$ in a natural way to  a continuous map $F: \psi (T_l) \to L$, for all $1\leq l \leq n$. Note that this extension is  also an affine map. So, we can extend this map to a diffeomorphism (that we call $F$) $C^{\infty}$ of $S^1$.

 The map $F$ defined in this way, restrict to $C_{R}$,  is conjugated to $f_{S}$ restrict to $C_{S}$ (by $\psi$). So it has a minimal Cantor set, and this is a contradiction.
\hfill $\square$

%
%
%

\end{document}